\documentclass[10pt]{article}
\usepackage{makeidx}         % allows index generation
\usepackage[bottom]{footmisc}% places footnotes at page bottom
\usepackage{url}
\usepackage{wrapfig,lipsum}
\usepackage{siunitx}
\usepackage[right]{lineno}
\usepackage{authblk}
\usepackage{blindtext}
\usepackage{amsmath}
\usepackage{amssymb}
\usepackage{cite}
\usepackage{lineno}
\usepackage{microtype}
\DisableLigatures[f]{encoding = *, family = * }
\usepackage{changepage}
\usepackage[utf8x]{inputenc}
\usepackage{textcomp,marvosym}
\usepackage{nameref}
\usepackage{hyperref} 
\usepackage{microtype}
\usepackage{graphicx}
\usepackage[toc,page]{appendix}
\usepackage{setspace} 
\usepackage{amsthm}
\usepackage{authblk}
\usepackage[labelfont=bf,labelsep=period,justification=raggedright]{caption}
\usepackage[font=small,labelfont=bf,labelsep=space]{caption}
\makeatletter
\renewcommand{\@biblabel}[1]{\quad#1.}
\makeatother

\begin{document}
 \title{\textbf{ Black holes and topological surgery}} 
\author[1]{Stathis Antoniou}
\author[2]{Louis H.Kauffman}
\author[3]{Sofia Lambropoulou}
\affil[1]{\small{School of Applied Mathematical and Physical Sciences, National Technical University of Athens, Greece}, \textit{santoniou@math.ntua.gr}}
\affil[2]{\small{Department of Mathematics, Statistics, and Computer Science, University of Illinois at Chicago, Chicago, USA, Department of Mechanics and Mathematics\\ Novosibirsk State University\\
              Novosibirsk, Russia}, \newline \textit{kauffman@uic.edu}}
\affil[3]{\small{School of Applied Mathematical and Physical Sciences, National Technical University of Athens, Greece}, \textit{sofia@math.ntua.gr}}
\date{}
\maketitle
	
%%%%%%%%%%%%%%%%%%%%%%%%%%%%%%%%%%%%%%%%%%%%%%%%%%%%%%%%%%%%%%%%%%%%%%%%%%%%%%%%%%
\section*{Abstract} 
We directly connect topological changes that can occur in mathematical three-space via surgery, with black hole formation, the formation of wormholes and new generalizations of these phenomena. This work widens the bridge between topology and natural sciences and creates a new platform for exploring geometrical physics.

%%%%%%%%%%%%%%%%%%%%%%%%%%%%%%%%%%%%%%%%%%%%%%%%%%%%%%%%%%%%%%%%%%%%%%%%%%%%%%%%%%
\newpage
\section*{Essay}
The universe undergoes topological and geometrical changes at all scales. This essay goes to the foundations of these changes by offering a novel topological perspective. We describe the formation of black holes and wormholes via topological surgery. Surgery is a manifold-changing process that creates new spaces from known ones and has been used in the study and classification of manifolds. In the cosmic scales, surgery allows the formation of a black hole from the collapse of a knotted cosmic string, without ending up in a singular manifold. 
It further describes Einstein-Rosen bridges (wormholes) linking the two black holes through a singularity where the disconnected black holes collapse to each other, and the bridge is born topologically. The collapse of a cosmic string can be viewed as an orchestrated creation of bridges that is topologically equivalent to $3$-dimensional surgery. We present a rich family of $3$-manifolds that can occur and the possible implications of these constructions in quantum gravity and general relativity. 

\begin{wrapfigure}{R}{0.4\textwidth}
\centering
\includegraphics[width=0.35\textwidth]{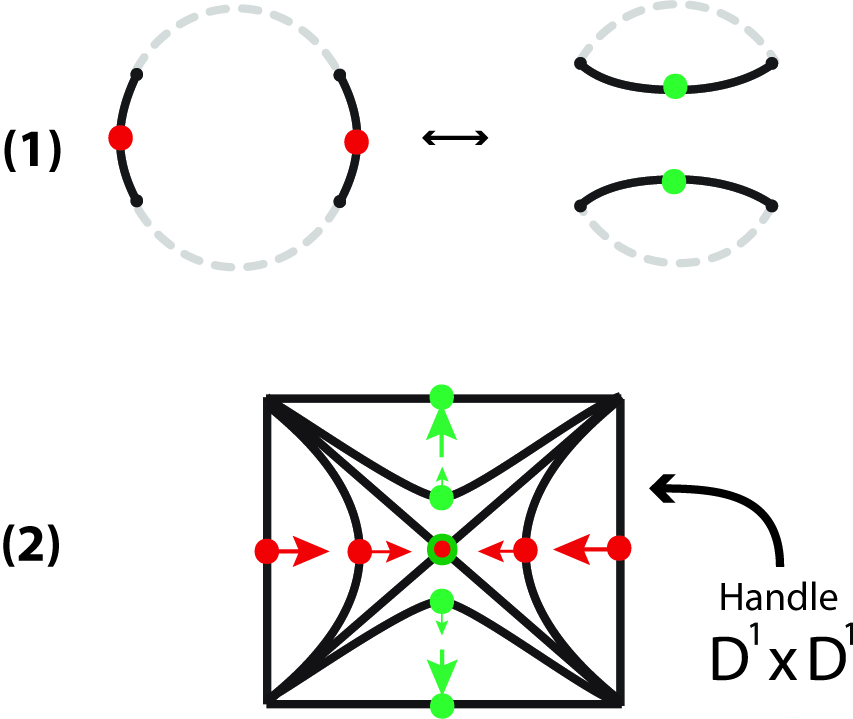}
\caption{\textbf{(1)} $1$-dimensional surgery \textbf{(2)} $2$-dimensional handle}
\label{SurgEx}
\end{wrapfigure}

The key surgery idea is to cut a thickened sphere from a manifold and to glue back another thickened sphere with the same boundary. The process of $1$-dimensional surgery causes the transformation of a circle into two circles, see Fig.~\ref{SurgEx}~(1). In the figure, the $1$-dimensional thickening of a $0$-sphere (represented by the two red points) is removed and replaced by the the thickening of another $0$-sphere (represented by two green points). This operation causes a global change which is induced by a local process. The local process happens within a handle which is of one dimension higher, see Fig.~\ref{SurgEx}~(2). This extra dimension leaves room for continuously passing from one boundary component of the handle to the other. This process is exhibited in various natural phenomena where segments are detached and rejoined, such as site-specific DNA recombination and magnetic reconnection, see~\cite{SS1}. The surgery approach describes both the occurring topological change and their temporal evolution. 

\begin{wrapfigure}{L}{0.5\textwidth}
\centering
\includegraphics[width=0.45\textwidth]{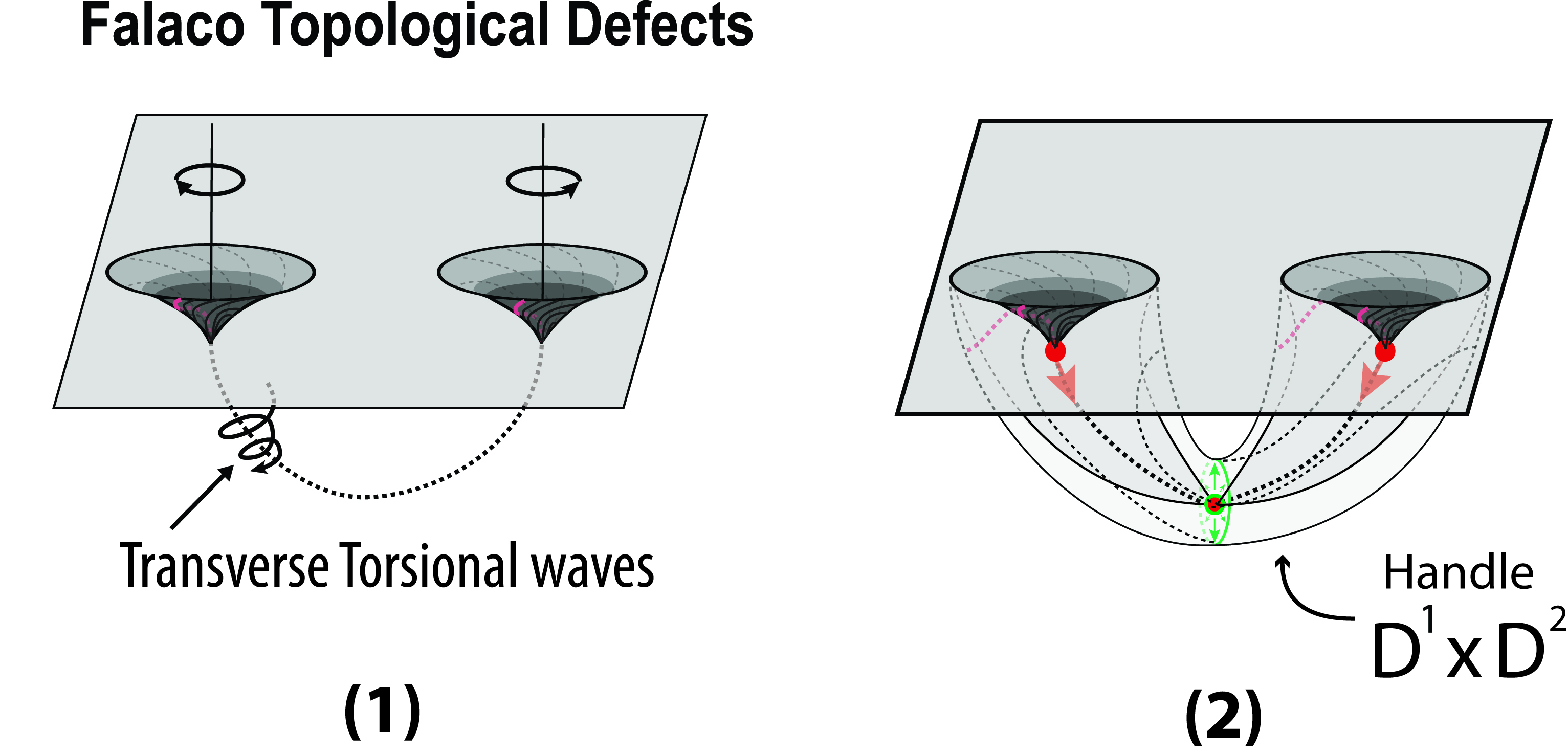}
\caption{\textbf{(1)} Falaco Solitons \textbf{(2)} $3$-dimensional handle}
\label{2D_Falacohomeo2}
\end{wrapfigure}

Nature is filled with $2$-dimensional surgeries too. One such phenomenon is the formation of Falaco solitons. Each Falaco soliton consists of a pair of contra-rotating identations in the water-air surface of a swimming pool, see Fig.~\ref{2D_Falacohomeo2}~(1) and \cite{Ki}. From the topological viewpoint the surgery consists in taking disk neighborhoods of two points (the identations in Fig.~\ref{2D_Falacohomeo2}~(1)) and joining them via a tube (which is a thickened circle $D^1 \times S^1$), see Fig.~\ref{2D_Falacohomeo2}~(2). Here the tube is the cylindrical vortex made from the propagation of the torsional waves around the singular thread. The $3$-dimensional handle containing all the $2$-dimensional temporal `slices' of this process is shown Fig.~\ref{2D_Falacohomeo2}~(2). 

One dimension higher, we have the two types of $3$-dimensional surgery, both of which require four dimensions to be visualized. The first one joins the spherical  neighborhoods of two points via a tube $D^1 \times S^2$ which is one dimension higher than the one shown in Fig.~\ref{2D_Falacohomeo2}~(2). If we consider that our initial manifold is the $3$-dimensional spatial section of the $4$-dimensional spacetime, this tube is what physicists call a wormhole. A connection between Falaco solitons and wormholes has been conjectured by R.M. Kiehn~\cite{KiehnSmall}. Our surgery description reinforces this connection. Moreover, this change of topology, which, according to J.A. Wheeler, results from quantum fluctuations at the Planck scale~\cite{Whee_QF}, can now also be viewed as a result of a ‘classical’ continuous topological change of $3$-space.

\begin{wrapfigure}{R}{0.5\textwidth}
\centering
\includegraphics[width=0.45\textwidth]{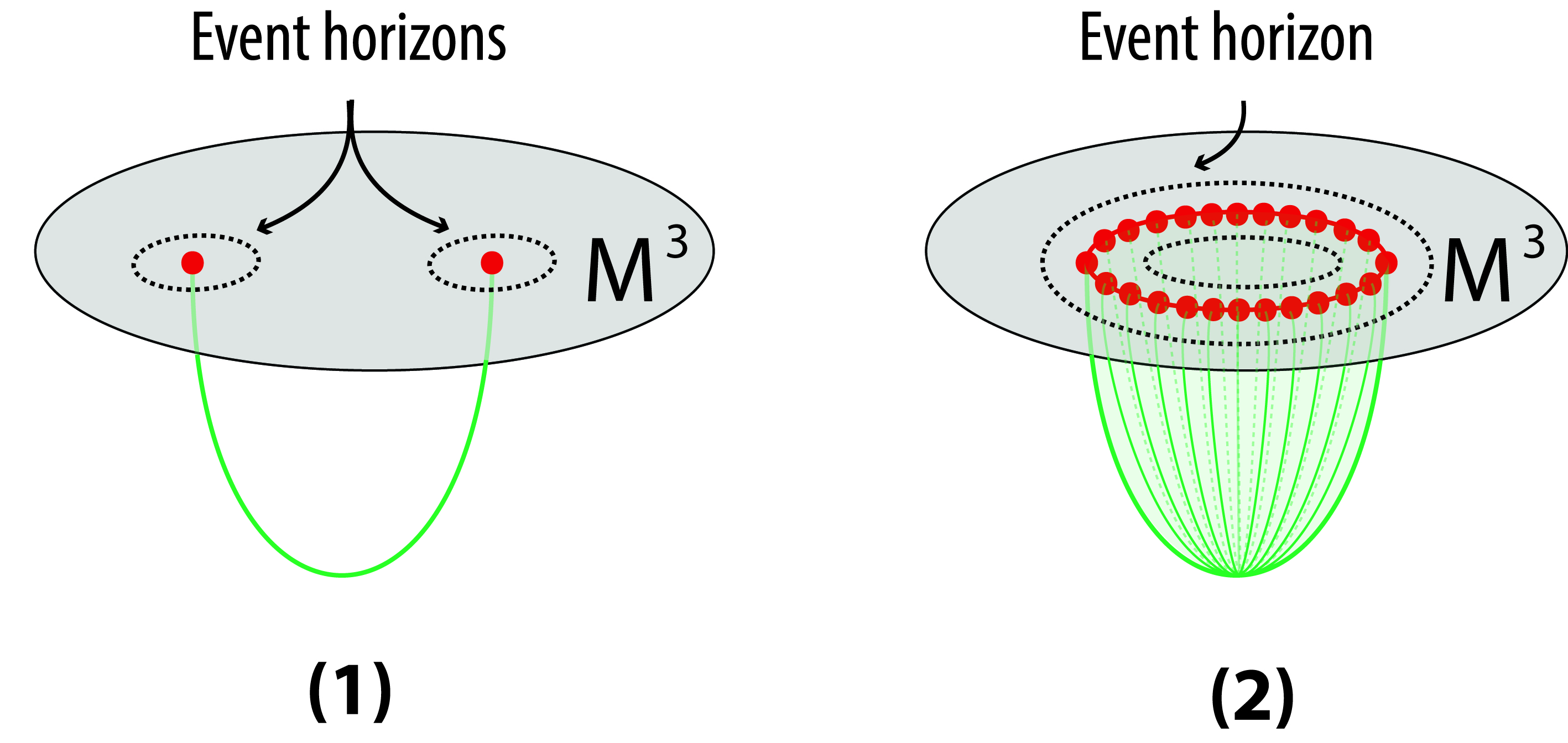}
\caption{\textbf{(1)} Pair of entangled black holes \textbf{(2)} String of entangled black holes}
\label{EnStrings}
\end{wrapfigure}

Let us now consider the $ER = EPR$ hypothesis of L. Susskind~\cite{ER_EPR}, which says that a wormhole is equivalent to the quantum entanglement of two concentrated masses that each forms its own black hole. Adding this hypothesis to our description, the two sites in space are the singularities of the two black holes, shown in red in Fig.~\ref{EnStrings}~(1), which will not collapse individually but will become the ends of the wormhole, shown in green in Fig.~\ref{EnStrings}~(1). We cannot visualize this process directly but it can be understood by considering that the green arc is the core $D^1$ of the higher dimensional handle $D^1 \times D^3$ for the wormhole. Note that an observer in our initial $3$-space $M^3$ would not be able to detect the topological change,  which occurs across the event horizons.

\begin{wrapfigure}{L}{0.5\textwidth}
\centering
\includegraphics[width=0.4\textwidth]{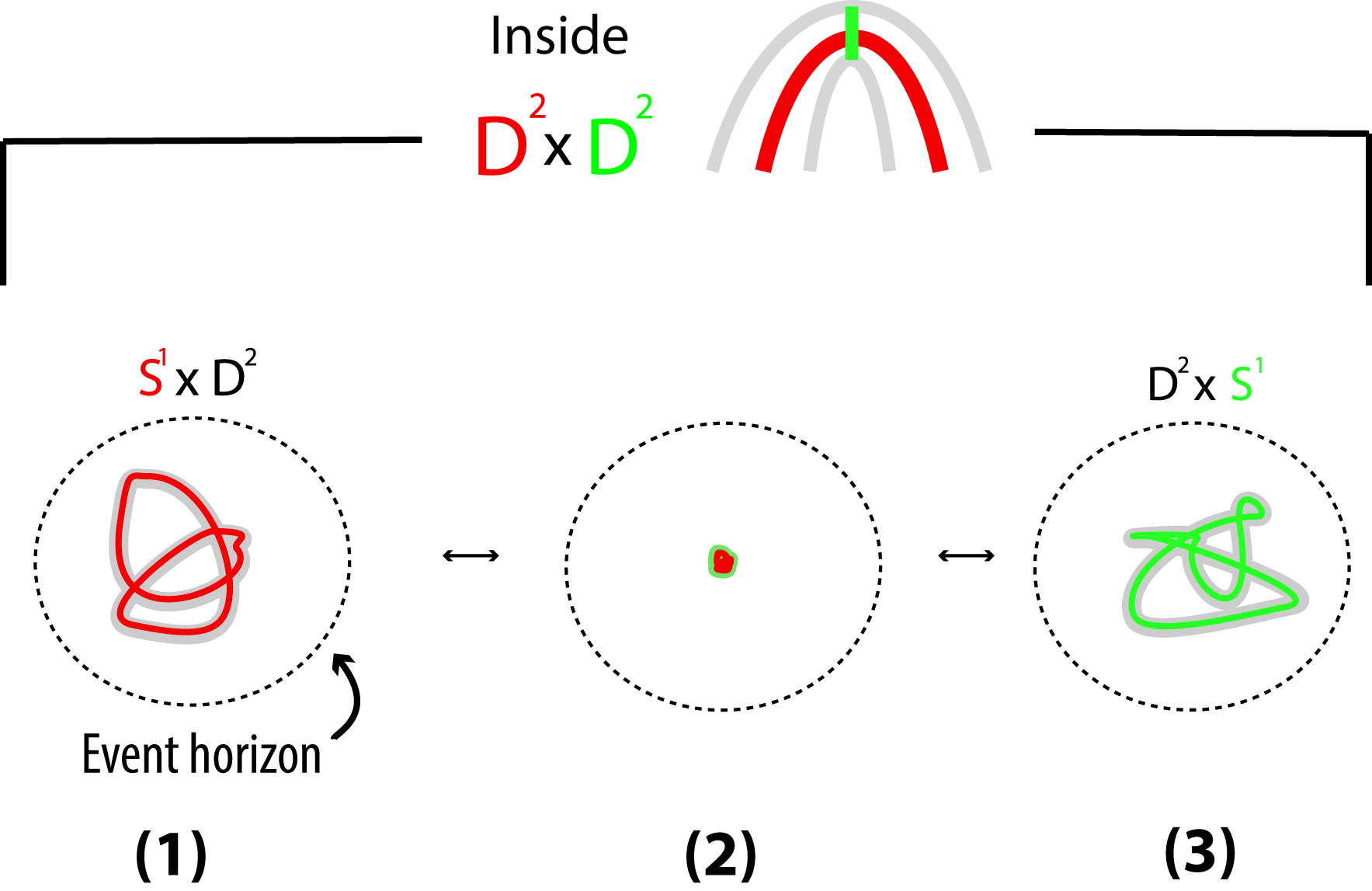}
\caption{3-dimensional surgery inside the event horizon}
\label{Fig35}
\end{wrapfigure}

The other type of $3$-dimensional surgery describes a more subtle topological change. It collapses a solid torus (which is a thickened circle) to a point and uncollapses another solid torus in such way that the meridians are glued to the longitudes and vice-versa. This type of surgery is also called `knot surgery' as the circle can be a knot. Knot surgery is an ideal candidate for describing black holes that are formed via the collapse of cosmic strings. This idea is based on~\cite{Hawk} where S.W. Hawking estimates that a fraction of cosmic string loops can collapse to a small size inside their Schwarzschild radius. As cosmic strings are hypothetical topological defects of small (but non-zero) diameter, a cosmic string loop can be considered as a knotted solid torus. As described in~\cite{Hawk}, the loop collapses to a point thus creating a black hole the center of which contains the singularity. At that point, the $3$-space becomes singular, see the passage from Fig.~\ref{Fig35}~(1) to~(2).    

Our surgery description says more~\cite{AKL}. According to it the process doesn't stop at the  singularity, but continues with the uncollapsing of another cosmic string loop from the singularity, see Fig.~\ref{Fig35}~(3). Thus, the creation of a cosmic string black hole is a $3$-dimensional surgery that continuously changes the initial $3$-space to another $3$-manifold. The process goes through the singular point of the black hole without having a singular manifold in the end. Instead, one ends up with a topologically new universe with a local topology change in the $3$-space, which happens within the event horizon.

This type of surgery is also related to the $ER = EPR$ hypothesis. Consider a cosmic string made of pairs of entangled concentrated masses. When each pair of masses collapse, they become connected by a wormhole as previously shown in Fig.~\ref{EnStrings}~(1). Given that all these pairs of masses have started on the same cosmic string, the distinct wormholes merge and the entire collection of wormhole cores (the green arcs in Fig.~\ref{EnStrings}~(1)) forms a 2-disc $D^2$, see Fig.~\ref{EnStrings}~(2), which is the core of the higher dimensional handle $D^2 \times D^2$ containing the temporal `slices' of the process. Our surgery description generalizes having a separate Einstein-Rosen bridge for each pair of black holes and amalgamates these bridges to form a new manifold in three dimensions. The effect of surgery is that, from any black hole location on the cosmic string to any other, there is a `bridge' through the new $3$-manifold. As this process joins the neighborhood of a circle instead of two points, one can rotate Fig.~\ref{EnStrings}~(1) to receive Fig.~\ref{EnStrings}~(2).

Another advantage of knot surgery is that it is able produce a great variety of $3$-manifolds. In fact, according a theorem by A.H. Wallace~\cite{Wal} and W. Lickorish~\cite{LickTh} knot surgery can create all closed, connected, orientable $3$-manifolds. One such $3$-manifold, which is of great interest to physicists, is the Poincaré dodecahedral space, that has been proposed as possible shape for the geometric universe~\cite{Weeks,Luminet,Levin}. This manifold can be obtained by doing knot surgery on the trefoil knot (with the right framing, see~\cite{PS}). Did the shape of the universe come about via the collapse of a trefoil cosmic string?!

\begin{wrapfigure}{R}{0.4\textwidth}
\centering
\includegraphics[width=0.35\textwidth]{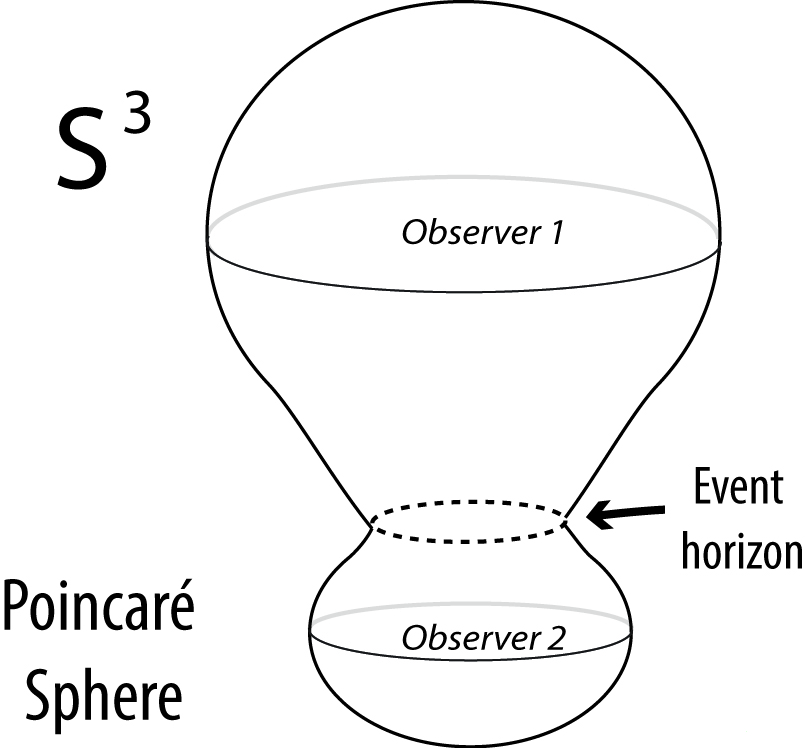}
\caption{Observer 1 and 2}
\label{BH_PoincaFlip}
\end{wrapfigure}

Suppose there are observers in an initial spherical universe $M^3=S^3$ containing a trefoil cosmic string. After surgery, a `mathematical' observer would be able to see the Poincaré dodecahedral space and detect the topology change.  However, a physical observer, who is subject to the restrictions of physical laws, would only see towards the event horizon in which the trefoil cosmic string has collapsed. Let us call this observer, \textit{Observer 1}, see Fig.~\ref{BH_PoincaFlip}. After surgery, \textit{Observer 1} would see the same universe $S^3$, the only change being the formation of the event horizon. On the other side of the event horizon, a new universe emerges in which new observers might evolve. Such an  observer, say \textit{Observer 2}, will see a Poincaré dodecahedral space and the event horizon from the other side, unaware that the original $S^3$ universe is behind it, see Fig.~\ref{BH_PoincaFlip}. Finding the Poincaré dodecahedral space (or some other non-trivial $3$-manifold) in our universe may indicate that we are observers that evolved inside the event horizon of a collapsed trefoil cosmic string (or some other cosmic string).

The surgery approach provides continuous paths to wormhole and cosmic string black hole formation. If one adds the $ER = EPR$ hypothesis, surgery also describes the entanglement of a pair or a string of black holes. Our topological perspective offers a process producing black holes and new non-singular $3$-manifolds from cosmic strings, binding entanglement and the connectivity of space with the rich structure of three- and four-dimensional manifolds.

%%%%%%%%%%%%%%%%%%%%%%%%%%%%%%%%%%%%%%%%%%%%%%%%%%%%%%%%%%%%%%%%%%%%%%%%%%%%%%%%%

\end{document}